\documentstyle[12pt,twoside]{article}
\font\twelvegtc=eufm10 scaled 1200

\font\ninegtc=eufm9
\font\sevengtc=eufm7

\newfam\gtcfam

\textfont\gtcfam=\twelvegtc
\scriptfont\gtcfam=\ninegtc
\scriptscriptfont\gtcfam=\sevengtc
\font\twelveBBB=msbm10 scaled 1200
\font\tenBBB=msbm10
\font\sevenBBB=msbm7
\newfam\BBBfam
\def\BBB{\fam\BBBfam\twelveBBB}
\textfont\BBBfam=\twelveBBB
\scriptfont\BBBfam=\tenBBB
\scriptscriptfont\BBBfam=\sevenBBB
\def\QQ{{\BBB Q}}       
   \def\CC{{\BBB C}}    \def\PP{{\BBB P}}

\def\hom{{\rm Hom}}
\def\ext{{\rm Ext}}

\newcounter{No}   \newcounter{SubNo}[No]    \newcounter{SubSubNo}[SubNo]
\renewcommand{\theNo}{\arabic{No}}
\renewcommand{\theSubNo}{\arabic{No}.\arabic{SubNo}}
\renewcommand{\theSubSubNo}{\arabic{No}.\arabic{SubNo}.\arabic{SubSubNo}}
\def\No#1{\refstepcounter{No}\vspace{3ex}\par
            \noindent {\large\bf\theNo.\hspace{2pt}#1.}\par}
\def\SubNopt#1{\refstepcounter{SubNo}\vspace{2ex}\par\noindent
                               {\bf\theSubNo.\hspace{2pt}#1.}\hspace{1ex}}
\def\SubNo#1{\refstepcounter{SubNo}\vspace{2ex}\par\noindent
                               {\bf\theSubNo.\hspace{2pt}#1}\hspace{1ex}}
\def\SubSubNo#1{\refstepcounter{SubSubNo}\vspace{1ex}\par\noindent
                               {\bf\theSubSubNo.}\hspace{2pt}#1.\hspace{1ex}}
\def\EF{\endgroup\par\vspace{1ex}\par}
\def\Pr{\refstepcounter{SubNo}\vspace{2ex}\par\noindent
       {\bf\theSubNo.}\hspace{2pt}PROPOSITION.\hspace{1ex}\begingroup\sl}
\def\Cl{\refstepcounter{SubNo}\vspace{2ex}\par\noindent
       {\bf\theSubNo.}\hspace{2pt}COROLLARY.\hspace{1ex}\begingroup\sl}
\def\proof{\par\noindent{\sc Proof.}\hspace{1ex}}
\def\Lm{\SubSubNo{LEMMA}\begingroup\sl}
\def\mst{{$\mu$-stable}}
\def\mss{{$\mu$-semistable}}
\def\gss{{$\gamma$-semistable}}
\def\gst{{$\gamma$-stable}}
\def\tE{\tilde E}
\def\tF{\tilde F}
\def\tx{\tilde x}
\title{On the structure of semistable rigid sheaves
on algebraic surfaces}
\author{B.V. Karpov \thanks{The work on this paper was supported by
the International Science Foundation grant MKU 300 and by the
Russian Foundation of Fundamental Research grant {\it N} 95-01-00840}}
\date{November 1995}
\begin{document}
\maketitle
\begin{abstract}
Let $S$ be a smooth projective surface, $K$ be the canonical class of $S$
and $H$ be an ample divisor such that $H\cdot K<0$. In this paper we
prove that for any rigid ($\ext^1(F,F)=0$) semistable sheaf $F$ in
the sense of Mumford--Takemoto stability w.r.t. $H$ there exists an
exceptional collection $(E_1,\dots ,E_n)$ of sheaves on $S$
such that $F$ can be constructed from $\{ E_i\}$ by a finite number of
extensions.
\end{abstract}
%
%
\markboth{B.V. Karpov}{Semistable rigid sheaves on surfaces}
\renewcommand{\theenumi}{(\roman{enumi})}

\par\vspace{2ex}
\noindent{\large\bf Introduction.}\par{}

\vspace{2ex}

Let $S$ be a smooth projective surface over  $\CC$ with the canonical
class $K$. Suppose $H$ is an ample divisor such that $H\cdot K<0$.

For any coherent sheaf $F$ on $S$ with $r(F)>0$ put by definition:
$$\mu(F)=\frac{c_1(F)\cdot H}{r(F)}\in\QQ\
\mbox{-- the Mumford--Takemoto slope of $F$ w.r.t. $H$;}
$$
$$\gamma(F)(n)=\frac{\chi(F(nH))}{r(F)}\in\QQ [n]
\mbox{ -- the Gieseker slope of $F$ w.r.t. $H$.}
$$
A coherent sheaf $F$ without torsion is $\mu$-{\it semistable}
({\it stable}) if for any subsheaf $E$ of $F$ such that
$0<r(E)<r(F)$ the following inequality holds:
$$\mu(E)\le\mu(F)\qquad (\mu(E)<\mu(F))\ .
$$
A torsion free sheaf $F$ is $\gamma$-{\it semistable}
({\it stable}) if for any subsheaf $E$ of $F$ with
$0<r(E)<r(F)$ the inequality
$$\gamma(E)(n)\le\gamma(F)(n)\qquad (\gamma(E)(n)<\gamma(F)(n))
$$
holds for $n\gg 0$.

A sheaf $F$ is {\it rigid} if $\ext^1(F,F)=0$ and
{\it superrigid} if $\ext^i(F,F)=0$ for $i>0$. An {\it exceptional} sheaf
is a superrigid sheaf $E$ such that $\hom(E,E)\cong\CC$.

An ordered collection of sheaves $(E_1,\dots,E_n)$ is called
{\it exceptional} if all $E_i$ are exceptional and
$\ext^i(E_k,E_j)=0,\ \forall i$ for $k>j$.

The aim of this paper is to prove the following result.

\vspace{2ex}

THEOREM. {\sl Any rigid {\mss} sheaf $F$ has a filtration
$$0=F_{n+1}\subset F_n\subset\cdots\subset F_2\subset F_1=F$$
such that all the quotients $F_i/F_{i+1}$ are isomorpfic to
$E_i\oplus\cdots\oplus E_i$ for some exceptional {\gst} sheaves $E_i$
and the collection
$(E_1,\dots,E_n)$ is exceptional.}

\vspace{2ex}

Results on the structure of rigid and superrigid sheaves are closely
connected with study of exceptional collections of sheaves over
del Pezzo surfaces (see [KuOr]) and on surfaces with nef anticanonical
class (see [Kul]). The theorem formulated above is a generalization of
[Kul,~2.4.1,~3.].

Actually, the existence of an ample divisor having negative
intersection with the canonical class imposes a restriction on the
surface. However, the main result of this paper gives a new information,
in particular, for rational ruled surfaces
$\PP({\cal O}_{\PP_1}\oplus {\cal O}_{\PP_1}(m))$, $m\ge 2$ and for
appropriate blowings up of them.

Observe that for $S=\PP({\cal O}_{\PP_1}\oplus {\cal O}_{\PP_1}(m))$
we have $n\le 3$ because the maximal possible length of an exceptional
collection is $\dim K_0(S)=4$ and the images of $E_i$ in $K_0(S)$
belong to the linear subspace $c_1\cdot H=r\mu(F)$.

The author is grateful to S.~A.~Kuleshov for stimulating discussions and
to A.~N.~Rudakov and A.~L.~Gorodentsev for attension to this
work.
\No{Stability}
\SubNo{} We begin with some notation. We write that $\gamma (E)>
\gamma (F)$ if $\gamma (E)(n)>\gamma (F)(n)$ for $n\gg 0$.
By the Riemann--Roch theorem,
$$\gamma (F)(n)=\frac{1}{2}n^2H^2+n\Bigl(\mu(F)+\frac{1}{2}H\cdot K
\Bigr)+\frac{\chi(F)}{r(F)}\quad .
$$
It is readily seen that
$$\arraycolsep=0.2em
\gamma (E)>\gamma (F)\quad\Longleftrightarrow\quad
\mu (E)>\mu (F)\:\mbox{ or }
\left\{ \begin{array}{ccl}
\mu (E)&=&\mu(F)\\[\medskipamount]
\frac{\textstyle\chi(E)}{\textstyle r(E)}&>&
\frac{\textstyle\chi(F)}{\textstyle r(F)}\quad .
\end{array}\right.
$$
This implies that any {\mst}  sheaf is also {\gst}  and any
{\gss}  sheaf is also {\mss}.

\SubNo{} We shall use the following properties of (semi)stable
sheaves.
\begin{description}
\item[\quad 1.] A torsion free sheaf $F$ is {\gss} ({\gst}) iff
$\gamma(F)\le\gamma(G)\ (\gamma(F)<\gamma(G))$
for any torsion free quotient sheaf $G$ of $F$.
\item[\quad 2.] If $E$ and $F$ are {\gss} and $\gamma (E)>\gamma (F)$, then
$\hom(E,F)=0$.
\item[\quad 3.] Suppose $0\longrightarrow F'\longrightarrow F
\longrightarrow G\longrightarrow 0$ is an exact triple of sheaves, where
$F$ and $G$ are {\gss} and $\gamma (F)=\gamma (G)$. Then
$\gamma (F')=\gamma (F)$ and $F'$ is {\gss}.
\item[\quad 4.] Let $E_1$,$E_2$ be {\gst} sheaves,
$\gamma (E_1)=\gamma (E_2)$, and $\hom(E_1,E_2)\ne 0$; then
$E_1\cong E_2$ and $\hom(E_1,E_2)\cong\CC$.
\item[\quad 5.] The first three statements are true for {\mss}
sheaves.
\end{description}

These properties are well known, see, for example, [Kul].
Besides, for (semi)\-stable sheaves on $\PP_n$ these statements
are proved in [OSS]. It is not hard to see that the proof in
our case is similar.

\Cl\label{e2}\ \ Suppose $E_1$ and $E_2$ are {\mss} sheaves,
${\mu (E_1)\ge\mu (E_2)}$; then $\ext^2(E_2,E_1)=0$.
\EF
\proof By Serre duality, it follows that
$\ext^2(E_2,E_1)\cong\hom^{\ast}(E_1,E_2\otimes K)$.
On the other hand, we have $c_1(E_2\otimes K)=c_1(E_2)+r(E_2)K$.
Consequently, we get
$\mu(E_2\otimes K)=\mu (E_2)+K\cdot H<\mu(E_2)\le\mu(E_1)$.
To conclude the proof, it remains to use the stability properties.
\No{Properties of f\/iltrations}
\SubNo{} Let us introduce some notation.
Let F be a sheaf. By  $xF$ we denote the sum
$\underbrace{F\oplus\cdots\oplus F}_{\mbox{$x$ times}}$.
Suppose
\begin{equation}\label{F}
0=F_{n+1}\subset F_n\subset\cdots\subset F_2\subset F_1=F
\end{equation}
is a filtration of $F$.  Denote by $Gr(F)$ the
ordered collection of quotients
$$Gr(F)=(F_n/F_{n+1},F_{n-1}/F_n,\dots ,F_1/F_2)\: .
$$

Sometimes we write that $F$ has a filtration with
$Gr(F)=(G_n,G_{n-1},\dots ,G_1)$. This means that there is a
filtration (\ref{F}) such that $G_i=F_i/F_{i+1},\
i=1,\dots ,n$. Note that $G_n=F_n$.

We shall use the spectral sequence
\begin{equation}\label{sp}
E_1^{p,q}=\bigoplus_{i}\ext^{p+q}(G_i,G_{i+p})
\end{equation}
which converges to $\ext^{p+q}(F,F)$.

Evidently, any filtration with two quotients is nothing
more than exact triple
$$0\longrightarrow G_2\longrightarrow F
\longrightarrow G_1\longrightarrow 0\: .$$

\vspace{-1ex}

\Lm {\rm (Mukai)} If $\hom(G_2,G_1)=\ext^2(G_1,G_2)=0$ in the
above-mentioned exact triple,  then
$\displaystyle\dim\ext^1(F,F)\ge\dim\ext^1(G_1,G_1)+
\dim\ext^1(G_2,G_2)$.
\EF
The proof is found in [Mu] or in [KuOr]. Besides, the
reader will easily prove this fact via the spectral sequence.
\Pr \label{join}
Let $F$ be a sheaf; then the following conditions are equivalent:
\begin{description}
\item[{\rm (i)}] there exists a filtration with
$Gr(F)=(G_n,\dots ,G_{i+1},G_i,\dots ,G_1)$;
\item[{\rm (ii)}] there exists a filtration with
$Gr(F)=(G_n,\dots ,G_{i+2},\tilde G,G_{i-1},\dots ,G_1)$,
where $\tilde G$ has a filtration with $Gr(\tilde G)=(G_{i+1},G_i)$.
\end{description}
\EF
\proof (i) $\Rightarrow$ (ii). Let
$\displaystyle 0\subset F_n\subset\cdots\subset F_{i+2}
\subset F_{i+1}\subset F_i\subset\cdots\subset F$
be the filtration with $Gr(F)=(G_n,\dots ,G_{i+1},G_i,\dots ,G_1)$;
then this filtration with omitted $F_{i+1}$ has the quotients
$Gr(F)=(G_n,\dots ,G_{i+2},\tilde
G,G_{i-1},\dots ,G_1)$. Here $\tilde G=F_i/F_{i+2}$. From the
commutative diagram
$$\arraycolsep=0.2em
\begin{array}{ccccccccc}
   & & 0 & & 0\\[\smallskipamount]
   & &\downarrow & &\downarrow \\[\smallskipamount]
   & &F_{i+2}    &= &F_{i+2}  \\[\smallskipamount]
   & &\downarrow & &\downarrow \\[\smallskipamount]
0&\longrightarrow &F_{i+1} &\longrightarrow &F_i &\longrightarrow
&G_i &\longrightarrow &0\\[\smallskipamount]
   & &\downarrow & &\downarrow & & \Vert\\[\smallskipamount]
0&\longrightarrow &G_{i+1} &\longrightarrow &\tilde G &\longrightarrow
&G_i &\longrightarrow &0\\[\smallskipamount]
   & &\downarrow & &\downarrow \\[\smallskipamount]
   & & 0 & & 0\\
\end{array}$$
it follows that $Gr(\tilde G)=(G_{i+1},G_i)$.

\noindent (ii)$\Rightarrow$ (i). Actually, we have only the middle column
and the lower row of the above-mentioned diagram. Clearly, one can reconstruct
the rest. This completes the proof.

\Pr \label{join2} Let
$\displaystyle 0\longrightarrow E\longrightarrow F
\stackrel{\alpha}{\longrightarrow}\tF\longrightarrow 0$
be an exact triple; then the following conditions
are equvivalent:
\begin{description}
\item[{\rm (i)}] $\tF$ has a filtration with
$Gr(\tF)=(G_n,\dots ,G_1)$;
\item[{\rm (ii)}] $F$ has a filtration with
$Gr(F)=(E,G_n,\dots ,G_1)$.
\end{description}
\EF
\proof We use induction on $n$. For $n=1$, there is nothing to
prove. Consider $n>1$.

\vspace{0.5ex}

\noindent (i) $\Rightarrow$ (ii). Suppose $\displaystyle
0\subset \tF_n\subset\cdots\subset \tF_2\subset \tF_1=\tF$
is the filtration with $Gr(\tF)=(G_n,\dots ,G_1)$; then $\tF_2$ has
the filtration with $Gr(\tF_2)=(G_n,\dots ,G_2)$. Consider the
epimorphism $\displaystyle
F\stackrel{\alpha}{\rightarrow}\tF\rightarrow G_1$ and its kernel
$F_2$. It is readily seen that there is the following commutative
diagram:
\begin{equation}\label{dia}
\arraycolsep=0.2em
\begin{array}{ccccccccc}
   & & 0 & & 0\\[\smallskipamount]
   & &\downarrow & &\downarrow \\[\smallskipamount]
   & &E    &= &E  \\[\smallskipamount]
   & &\downarrow & &\downarrow \\[\smallskipamount]
0&\longrightarrow &F_2 &\longrightarrow &F &\longrightarrow
&G_1 &\longrightarrow &0\\[\smallskipamount]
   & &\downarrow & &\downarrow & & \Vert\\[\smallskipamount]
0&\longrightarrow &\tF_2 &\longrightarrow &\tilde F &\longrightarrow
&G_1 &\longrightarrow &0\\[\smallskipamount]
   & &\downarrow & &\downarrow \\[\smallskipamount]
   & & 0 & & 0\\
\end{array}
\end{equation}
By the induction hypothesis, $F_2$ has a filtration with
$Gr(F_2)=(E,G_n,\dots ,G_2)$. From the middle row of the diagram
it follows that $F$ has the required filtration.

\vspace{0.5ex}

\noindent (ii) $\Rightarrow$ (i). Let
$0\subset E\subset F_n\subset\cdots\subset F_2\subset F$
be the filtration with
$Gr(F)=(E,G_n,\dots ,G_1)$;
then we have the two upper rows and the middle column of
(\ref{dia}). Clearly, the rest part of the diagram
can be reconstructed. By the induction
hypothesis, there is a filtration of $\tF_2$ with
$Gr(\tF_2)=(G_n,\dots ,G_2)$. Thus, the lower row of
diagram (\ref{dia}) implies that $\tF$ has the required filtration.
This concludes the proof.

\vspace{0.5ex}

Combining the two previous statements, we get the following.

\Pr\label{FF} A sheaf $F$ has a filtration with $Gr(F)=(G_n,\dots ,G_1)$
and there are filtrations $Gr(G_i)=(E^{(i)}_{m_i},\dots ,E^{(i)}_1)$
for $i=1, \dots ,n$ if and only if $F$ has a filtration with
$$Gr(F)=\Bigl(E^{(n)}_{m_n},\dots ,E^{(n)}_1\, ,\,E^{(n-1)}_{m_{n-1}},\dots ,
E^{(n-1)}_1,\:\dots\: ,\: E^{(1)}_{m_1},\dots ,E^{(1)}_1\Bigr)\ \: .
$$
\EF

\vspace{-0.5ex}

The proof is left to the reader.
\No{The case of {\gss} sheaf}
\SubNo{}\label{fgs} Suppose $F$ is {\gss} sheaf.
In the work [Kul, 1.3.7] it is proved that $F$ has a filtration
$$0=F_{n+1}\subset F_n\subset\cdots\subset F_2\subset F_1=F
$$
with {\gss} quotients
$$Gr(F)=(G_n,G_{n-1},\dots,G_1)\ ,
$$
where each $G_i$ has a filtration with stable quotients isomorphic
to each other:
$$Gr(G_i)=(E_i,E_i,\dots,E_i)\: ,\ \gamma(E_i)=\gamma(G_i)=
\gamma(F)\: .$$
Moreover, this filtration can be constructed such that
$$\hom(F_i,G_{i-1})=\hom(G_i,G_{i-1})=\hom(G_{i-1},G_i)=0\  .
$$

\vspace{-1.5ex}

\Pr\label{fgss} Let $F$ be rigid; then the following statements hold:
\par\vspace{1ex}
\noindent {\rm a)} all $F_i$ and $G_i$ are rigid;
\par\vspace{1ex}
\noindent {\rm b)} $\ext^1(G_n,G_1)=0$;
\par\vspace{1ex}
\noindent {\rm c)} $E_i$ are exceptional and $G_i=E_i\oplus\cdots
\oplus E_i$;
\par\vspace{1ex}
\noindent {\rm d)} $\ext^1(E_n,E_i)=0,\quad 1\le i\le n-1$.
\EF
\proof a) The sheaves $F_2$ and $G_1$ are semistable and
$\gamma(F_2)=\gamma(G_1)$. Using \ref{e2}, we get $\ext^2(G_1,F_2)=0$.
{}From Mukai Lemma it follows that $F_2$ and $G_1$ are rigid.
Taking $F_2$ in place of $F$, we similarly obtain that
$F_3$ and $G_2$ are rigid, etc.

\noindent b) Consider the spectral sequence associated with the
filtration. This sequence has
$E_1^{p,q}=\bigoplus_{i}\ext^{p+q}(G_i,G_{i+p}) $ and converges to
$\ext^{p+q}(F,F)$. Using \ref{e2}, we get $\ext^2(G_i,G_j)=0,\
\forall i,j$. Therefore, the first term of the spectral sequence
has the form:
$$\arraycolsep=0.2em
\begin{array}{cccccc}
       \  \  0     \\[\smallskipamount]
 \ext^1(G_n,G_1) &  0   &               \\[\smallskipamount]
 \hom(G_n,G_1)   & \ast &              &  0   &        &        \\
    & \ast &  & \ast &\quad\ddots &        \\
    &      &  & \ast &\quad\ddots &\raisebox{1em}[0.1em][0.2em]{$p+q=2$}\\
    &      &  &      &\quad\ddots &\raisebox{1.1em}[0.1em][0.2em]{$p+q=1$}\\
 &\cdot&\stackrel{d}{\longrightarrow}& \cdot & &
\raisebox{1em}[0.1em][0.2em]{$p+q=0$}\\
\end{array}$$

The sheaf $F$ is rigid whence the spectral sequence converges to $0$
on the diagonal $p+q=1$. Thus we have $\ext^1(G_n,G_1)=0$.

\noindent c) It follows from \ref{e2} that $\ext^2(E_i,E_i)=0$.
Taking $G_i$ in place of $F$ in b), we get $\ext^1(E_i,E_i)=0$.
This implies that $G_i=E_i\oplus\cdots\oplus E_i$. Finally,
$E_i$ is {\gst} whence $\hom(E_i,E_i)\cong\CC$.

The statement d) follows from c) and b) applied to $F_i$.
This concludes the proof.

\SubNopt{Key Lemma} Let $F$ be {\gss} rigid sheaf. Consider the
filtration
$$0=F_{n+1}\subset F_n\subset\cdots\subset F_2\subset F_1=F
$$
described in \ref{fgs}. By the previous proposition, this filtration has
the quotients
$$Gr(F)=(x_nE_n,\dots ,x_1E_1)\ ,\quad \gamma(E_i)=\gamma(F)\ ,
$$
where the sheaves $E_i$ are exceptional, {\gst}, and
$E_i\not\cong E_{i-1},\ i=2,\dots,n$. In addition,
$\hom(F_i,E_{i-1})=0$.

Suppose that there are some $i\in \{1,\dots,n-1\}$ such that
$\hom(E_n,E_i)\ne 0$.

Take the greatest $j$ such that $\hom(E_n,E_j)\ne 0$. From
the stability properties it follows that $E_j\cong E_n$.
Using \ref{fgss} d) and
applying the functor $\hom(E_n,\cdot)$ to the exact triples
$$\arraycolsep=0.2em
\begin{array}{ccccccccc}
 0 &\longrightarrow& x_n E_n      &\longrightarrow&    F_{n-1}   &
     \longrightarrow& x_{n-1}E_{n-1} &\longrightarrow&     0
\\[\smallskipamount]
 0 &\longrightarrow& F_{n-1}         &\longrightarrow& F_{n-2}     &
     \longrightarrow& x_{n-2}E_{n-2} &\longrightarrow&     0
\\[\smallskipamount]
   &&\vdots&&\vdots&&\vdots\\[\smallskipamount]
 0 &\longrightarrow&  F_{j+1}      &\longrightarrow& F_j          &
     \longrightarrow& x_j E_n      &\longrightarrow&     0
\\
\end{array}$$
we get: $\ext^1(E_n,F_{n-1})=\ext^1(E_n,F_{n-2})=\cdots=
         \ext^1(E_n,F_j)=0$.

Therefore, the last extension is trivial and $F_j=F_{j+1}\oplus x_jE_n$.
Whence, we have a new filtration of $F$
\begin{equation}\label{newfi}
\begin{array}{c}
0\subset (x_n+x_j)E_n\subset F_{n-1}\oplus x_jE_n\subset\cdots\subset
  F_{j+2}\oplus x_jE_n\subset\\[\medskipamount]
\subset F_j\subset F_{j-1}\subset\cdots\subset F
\end{array}\end{equation}
with the quotients
$$Gr(F)=\Bigl((x_n+x_j)E_n,x_{n-1}E_{n-1},\,\dots\, ,x_{j+1}E_{j+1},
x_{j-1}E_{j-1},\,\dots\, ,x_1E_1\Bigr)
$$
If $E_{j+1}\cong E_{j-1}$, then by \ref{join} we can join the
quotients $x_{j+1}E_{j+1}$ and $x_{j-1}E_{j-1}$ into the sum
$(x_{j+1}+x_{j-1})E_{j-1}$. The corresponding filtration is
(\ref{newfi}) with omitted $F_j$. Note that
$\hom(F_i\oplus x_jE_n,E_{i-1})=0$ for $i=j+2,\dots,n$.

Iterating this procedure, we obtain the following.

\vspace{1ex}\par

LEMMA.\hspace{1ex}{\sl
 There exists a filtration
$$0\subset{\tilde x}_n E_n={\tilde F}_{m+1}\subset{\tilde F}_m
\cdots\subset{\tilde F}_2\subset{\tilde F}_1=F
$$
with the quotients ${\tilde {Gr}}(F)=({\tilde x}_nE_n,
{\tilde x}_m{\tilde E}_m,\dots ,{\tilde x}_1{\tilde E}_1)$ such that
the ordered collection $({\tilde E}_1,\dots ,{\tilde E}_m)$ is a
subcollection of $(E_1,\dots ,E_{n-1})$ and
$$\hom(E_n,{\tilde E}_i)=0\ ,\quad i=1,\dots ,m\: .
$$
Furthermore, ${\tilde E}_m=E_{n-1}$,$\hom({\tilde F}_i,{\tilde E}_{i-1})=0$,
and $\tE_i\not\cong\tE_{i-1}$,$2\le i\le m$.
}

\No{Proof of the theorem}

\vspace{1ex}

The proof is in two steps.

\vspace{1ex}

\noindent STEP 1. {\bf $F$ is {\gss}}. We start with the filtration
described in the previous section:
$$0=F_{n+1}\subset F_n\subset\cdots\subset F_2\subset F_1=F\ .
$$
The quotients are:
$$Gr(F)=(x_nE_n,\dots ,x_1E_1)\ ,\quad \gamma(E_i)=\gamma(F)\ ,
$$
where the sheaves $E_i$ are exceptional and {\gst},
$\ext^1(E_n,E_i)=0$, $E_i\not\cong E_{i-1}$, and
$\hom(F_i,E_{i-1})=0$.
By the Key Lemma, we can construct this filtration such that
$\hom(E_n,E_i)=0$ for $i=1,\dots ,n-1$.

The sheaf $\tF =F/x_nE_n$ is {\gss} and has the filtration
\begin{equation}\label{tF}
0=F_n/x_nE_n\subset F_{n-1}/x_nE_n\subset\cdots\subset
F_2/x_nE_n\subset F_1/x_nE_n=\tF
\end{equation}
with {\gss} quotients
$$Gr(\tF)=(x_{n-1}E_{n-1},\dots ,x_1E_1)\ .
$$
Therefore, we obtain $\hom(E_n,\tF)=0$. Using \ref{e2},
we have $\ext^2(\tF,E_n)=0$. Applying Mukai Lemma to the
triple
$$0\longrightarrow x_nE_n\longrightarrow F\longrightarrow\tF
\longrightarrow 0\ ,$$
we conclude that $\tF$ is rigid.

Applying the functor $\hom(\cdot ,E_{i-1})$ to the triple
$$0\longrightarrow x_nE_n\longrightarrow F_i\longrightarrow
F_i/x_nE_n\longrightarrow 0\ ,$$
we get $\hom(F_i/{x_n E_n},E_{i-1})=0$. Hence the filtration
(\ref{tF}) of the sheaf $\tF$ is as in \ref{fgs}. From \ref{fgss} d) it
follows that
$$\ext^1(E_{n-1},E_i)=0,\ i=1,\dots ,n-2\: .
$$
Using the Key Lemma, we can construct a filtration of $\tF$
with the quotients
$$Gr(\tF)=(\tx_{n-1}E_{n-1},\tx_m\tE_m,\dots ,\tx_1\tE_1)
$$
such that the collection $(\tE_1,\dots ,\tE_m)$ is a subcollection
of $(E_1,\dots ,E_{n-2})$ (really, $\tE_m=E_{n-2}$) and
$\hom(E_{n-1},\tE_i)=0,\ i=1,\dots ,m$.
{}From \ref{join2} it follows that there exists a filtration of $F$
with the quotients
$$Gr(F)=(x_nE_n,\tx_{n-1}E_{n-1},\tx_m\tE_m,\dots ,\tx_1\tE_1)\ .
$$

Replacing $\tF$ by $\tF/\tx_{n-1}E_{n-1}$ in this reasoning and
iterating this procedure, we get a filtration of $F$ with
$$Gr(F)=(x_nE_n,\tx_{n-1}E_{n-1},\tx_{n-2}E_{n-2},x'_kE'_k,\dots ,
x'_1E'_1)\ ,
$$
where the collection $(E'_1,\dots ,E'_k,E_{n-2},E_{n-1},E_n)$ is
exceptional.

\vspace{2ex}

\noindent STEP 2. {\bf $F$ is {\mss} and rigid.} Let
$0\subset F_n\subset\cdots\subset F_2\subset F_1=F
$
be the Harder--Narasimhan filtration of $F$ w.r.t.$\:\gamma$ (see
[Kul]) with the quotients $Gr(F)=(G_n,\dots ,G_1)$. Here $G_i$ are
{\gss} and $\gamma (G_i)>\gamma (G_{i-1})$. The sheaf $G_n$ is a
subsheaf of $F$ and $G_1$ is a quotient sheaf of $F$ whence
$\mu (G_n)\le\mu (F)\le\mu (G_1)$. Therefore, we have
$\mu (G_i)=\mu (F)$ . Hence $\mu (F_i)=\mu (F)$ and $F_i$ are
semistable.

Let us show that $G_i$ are rigid. The sheaf $F_2$ has the filtration
with $Gr(F_2)=(G_n,\dots ,G_2)$. By the stability properties,
$\hom(G_i,G_1)=0,\ i=2,\dots ,n$. This implies that $\hom(F_2,G_1)=0$.
By \ref{e2}, we have $\ext^2(G_1,F_2)=0$. The Mukai Lemma implies
that $G_1$ and $F_2$ are rigid. Taking $F_3$ in place of $F_2$, we
similarly obtain that $G_2$ and $F_3$ are rigid. Continuing this
line of reasoning, we see that all $G_i$ are rigid.

By the first step, each $G_i$ has a filtration with
$$Gr(G_i)=(x^{(i)}_{m_i}E^{(i)}_{m_i},\dots ,x^{(i)}_1E^{(i)}_1)\: ,
\ i=1,\dots ,n\, ,$$
where $E^{(i)}_s$ are {\gst}, $\gamma (E^{(i)}_s)=\gamma (G_i),\
s=1,\dots ,m_i$, and the collection
$(E^{(i)}_1,\dots ,E^{(i)}_{m_i})$
is exceptional.
Proposition \ref{FF} implies that $F$ has a filtration with
\begin{equation}\label{final}
Gr(F)=\left( x^{(n)}_{m_n}E^{(n)}_{m_n},\dots ,x^{(n)}_1E^{(n)}_1,\:\dots\: ,
x^{(1)}_{m_1}E^{(1)}_{m_1},\dots ,x^{(1)}_1E^{(1)}_1\right)
\end{equation}
We claim that the collection
$\displaystyle
(E^{(1)}_1,\dots ,E^{(1)}_{m_1},\:\dots\: ,E^{(n)}_1,\dots ,E^{(n)}_{m_n})
$
is exceptional. By construction, $\gamma (E^{(j)}_t)>\gamma (E^{(i)}_s)$
for $j>i$, whence $\hom(E^{(j)}_t,E^{(i)}_s)=0$. On the other hand,
$\mu (E^{(j)}_t)=\mu (E^{(i)}_s)$. This implies that
$\ext^2(E^{(j)}_t,E^{(i)}_s)=0$. Finally, from the spectral sequence
(\ref{sp}) associated with the filtration (\ref{final}) it follows that
$\ext^1(E^{(j)}_t,E^{(i)}_s)=0$. This completes the proof.

\par\noindent
\begin{flushleft}
 {\large\bf References.}
\end{flushleft}
\vspace{1ex}
\begin{description}
\item{[KuOr]:}
S.A.Kuleshov, D.O.Orlov.
{\it Exceptional sheaves on Del Pezzo surfaces.}
Izv. Russ. Acad. Nauk Ser. Mat. Vol. {\bf 58} (1994), No.3, 59-93 (in Russian)
\item{[Kul]:}
S.A.Kuleshov.
{\it Exceptional and rigid sheaves on surfaces with anticanonical
class without base components.}
Preprint No.1 of Math. College of the Independent Moscow University, 1994 (in
Russian).
English version is available from e-print servise at
alg-geom@eprints.math.duke.edu, alg-geom/9511016.
\item{[Mu]:}
S.Mukai.
{\it On the Moduli Spaces of Bundles on K3 Surfaces, I.}
In: Vector Bundles ed. Atiyah et al, Oxford Univ. Press (1986), 67-83.
\item{[OSS]:}
C. Okonek, M. Schneider, H. Spindler.
{\it Vector Bundles on Complex Projective Spaces.}
Birkh\"auser, Boston, Basel, Stuttgart 1980.
\end{description}
\end{document}